\documentclass[11pt,reqno]{amsart}
\usepackage{amsfonts,amssymb,amsmath,amsthm,bbm}
\usepackage{enumitem}
\setlist[enumerate]{format=\normalfont}
\usepackage[colorlinks]{hyperref}

\usepackage{versions}

\includeversion{arXiv}
\excludeversion{arXiv-cut}

\usepackage[colorlinks]{hyperref}

\theoremstyle{plain}
\newtheorem{theorem}{Theorem}
\newtheorem*{proposition}{Proposition}
\newtheorem*{assumption}{Assumption}

\theoremstyle{remark}
\newtheorem*{remark}{Remark}
\newtheorem*{acknowledgements}{Acknowledgements}

\theoremstyle{definition}
\newtheorem*{definition}{Definition}

\usepackage{setspace}
\newcommand{\marginparstretch}{0.8}
\let\oldmarginpar\marginpar
\renewcommand\marginpar[1]{\-\oldmarginpar[\framebox{\setstretch{\marginparstretch}\begin{minipage}{\marginparwidth}{\raggedleft\scriptsize #1}\end{minipage}}]{\framebox{\setstretch{\marginparstretch}\begin{minipage}{\marginparwidth}{\raggedright\scriptsize #1}\end{minipage}}}}

\renewcommand{\baselinestretch}{1.2} 

\usepackage{tikz}
\usepackage{xcolor}
\usetikzlibrary{decorations.markings,arrows,calc,positioning} 

\tikzset{arrow/.style={decoration={markings,mark=at position 1 with %
    {\arrow[scale=1.7]{>}}},postaction={decorate}}}

\newcommand\twistLocal{\mathsf{T}_{\! E}} 
\newcommand\twistGlobal{\mathsf{T}_{\hspace{-0.1em}\mathcal{E}}}

\begin{document}

\begin{arXiv}\title[Applications of noncommutative deformations]{\phantom{text}\vspace{-1.7cm}Applications of \\ noncommutative deformations}\end{arXiv} 
\begin{arXiv-cut}\title[Applications of noncommutative deformations]{\phantom{text}\vspace{-1.2cm}Applications of \\ noncommutative deformations}\end{arXiv-cut} 
\author{W.\ Donovan}
\address{Kavli IPMU (WPI), UTIAS, University of Tokyo}
\email{will.donovan@ipmu.jp}
\begin{abstract}
For a general class of contractions of a variety $X$ to a base~$Y$\!, I discuss recent joint work with M.~Wemyss defining a noncommutative enhancement of the locus in $Y$\! over which the contraction is not an isomorphism, along with applications to the derived symmetries of $X$.  \begin{arXiv}This note is based on a talk given at the Kinosaki Symposium in~2016.\end{arXiv}
\end{abstract}

\thanks{The author is supported by World Premier International Research Center Initiative (WPI), MEXT, Japan, and JSPS KAKENHI Grant Number~JP16K17561.}
\maketitle
\tableofcontents

Derived symmetry groups of algebraic varieties extend classical symmetry groups to include contributions from symplectic geometry via homological mirror symmetry, and from birational geometry. In a recent joint paper~\cite{DW4}, M.~Wemyss and I construct, for a general class of birational contractions $f\colon X \to Y$\!,\, a sheaf of noncommutative algebras $\mathcal{A}$ on~$Y$ which induces a derived symmetry of~$X$ in an appropriate crepant setting. This short note explains key features of our results.

The sheaf $\mathcal{A}$ is supported on the locus of\, $Y$ over which $f$ is not an isomorphism. In previous joint work~\cite{DW1,DW3} we considered contractions of $3$-folds for which this locus is just a point. In this setting we studied an algebra of noncommutative deformations~$A$ which allowed new constructions of derived symmetries, and extended and unified known invariants of such contractions. I~begin by reviewing this, as $\mathcal{A}$ may be viewed as a sheafy version of the algebra~$A$.

\newpage
I also briefly discuss an example in which $f$ is a Springer resolution~(\S\ref{section springer}), and indicate recent work in which deformation algebras are used to recover the geometry of contractions~(\S\ref{section conjectures}).

\begin{acknowledgements}{}

I am grateful to the organisers and supporters of the Kinosaki Symposium for the opportunity to take part in the fine tradition of this meeting. My thanks also go to M.~Wemyss, as our joint work forms the subject of this note. The presentation of this work has benefited from comments from many people: in particular, I am grateful for recent conversations with 
A.~Bodzenta,
A.~Bondal, 
Z.~Hua,
Y.~Kawamata,
S.~Mehrotra,
T.~Logvinenko, 
D.~Piyaratne,
and Y.~Toda. 

\end{acknowledgements}

\subsubsection*{Conventions} I work over the ground field $\mathbb{C}$, though this assumption can be weakened. Varieties $X$ are assumed quasi-projective, with bounded derived category of coherent sheaves denoted by $D(X)$. The variety of hyperplanes in a vector space $V$ is denoted by $\mathbb{P}V$.

\section{{Deformation algebras for $3$-folds}}

The theorem below applies noncommutative deformations to study derived symmetries of $3$-folds. Given smooth $3$-folds $X$ and $X'$ related by a flop, Bridgeland \cite{BriFlop} constructs certain canonical Fourier--Mukai equivalences
\begin{equation*}\label{equation flops}\begin{tikzpicture}
\node (A) at (0,0) {$\phantom{.}D(X)$};
\node (B) at (3,0) {$D(X').$};
\node (C) at (1.5,0) {$\scriptstyle\sim$};
\draw[arrow,transform canvas={yshift=+0.8ex}] (A) -- node[above] {$\scriptstyle\mathsf{F}$} (B);
\draw[arrow,transform canvas={yshift=-0.6ex}] (B) -- node[below] {$\scriptstyle\mathsf{F'}$} (A);
\end{tikzpicture}
\vspace{-0.25cm}
\end{equation*}
These equivalences are not mutually inverse: the theorem explains this using deformations of curves on $X$.

Consider a $3$-fold $Y$\ with an isolated rational singular point~$p$, and a resolution  $f\colon X \to Y$ of this singularity, with one-dimensional exceptional locus with components~$C_i$ for $i=1, \dots, n$.

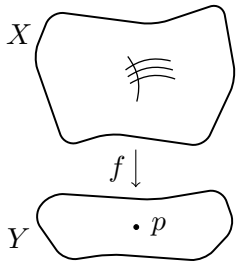
\begin{figure}[h]
\begin{center}
\begin{tikzpicture}
\node at (0,0) {\begin{tikzpicture}[scale=0.6]
\coordinate (T) at (1.9,2);
\coordinate (B) at (2.1,1);
\coordinate (L) at (1.8,1.4);
\coordinate (R) at (2.8,1.5);
\draw[line width=0.5pt] (T) to [bend left=25] (B);
\draw[line width=0.5pt] ($(L)+(0.15,0)$) to [bend left=25] ($(R)+(0.15,0)$);
\draw[line width=0.5pt] ($(L)+(0.1,0.15)$) to [bend left=25] ($(R)+(0.1,0.2)$);
\draw[line width=0.5pt] ($(L)+(0.05,0.3)$) to [bend left=25] ($(R)+(0.05,0.4)$);
\draw[rounded corners=5pt,line width=0.75pt] (0.5,0) -- (1.5,0.3)-- (3.6,0) -- (4.3,1.5)-- (4,3.2) -- (2.5,2.7) -- (0.2,3) -- (-0.2,2)-- cycle;
\node at (-0.5,2.5) {$X$};
\end{tikzpicture}};
\node at (0,-2) {\begin{tikzpicture}[scale=0.6]
\filldraw (2.1,0.75) circle (1.5pt);

\draw[rounded corners=5pt,line width=0.75pt] (0.5,0) -- (1.5,0.15)-- (3.6,0) -- (4.3,0.75)-- (4,1.6) -- (2.5,1.35) -- (0.2,1.5) -- (-0.2,1)-- cycle;
\node at (-0.5,0.5) {$Y$};
\node (label) at (2.6,0.75) {$p$};

\end{tikzpicture}};
\draw[->] (0.25,-1) -- node[color=black,left]{$f$} (0.25,-1.5);
\end{tikzpicture}
\end{center}
\caption{Contraction $f$ for Theorem~\ref{theorem flop flop}.}
\end{figure}

\begin{theorem}\label{theorem flop flop} \cite{DW1,DW2,DW3} Noting that the subvarieties $C_i$ of $X$\! are projective lines, we have that:
 
\begin{enumerate}
\item\label{theorem flop flop 1} there exists a $\mathbb{C}^n$-algebra $A$ which represents the functor of noncommutative deformations of the sheaves~$\mathcal{O}_{C_i}(-1)$ on $X$.
\end{enumerate}

\noindent Write $E$ for the corresponding universal sheaf on $X$. If the contraction~$f$ corresponds to a flop of $X$, then:

\begin{enumerate}[resume]
\item\label{theorem flop flop 2} there is a Fourier--Mukai autoequivalence $\twistLocal$ of~$D(X)$, fitting into a distinguished triangle of functors
\begin{equation*}\label{triangle} \mathbb{R}\!\operatorname{Hom}_X(E,-) \underset{A}{\overset{\mathbb{L}}{\otimes}} E \longrightarrow \operatorname{Id}_{D(X)} \longrightarrow \twistLocal \longrightarrow; \end{equation*}
\item\label{theorem flop flop 3} there is a natural isomorphism of functors
\[\twistLocal \cong \left( \mathsf{F}' \circ \mathsf{F} \right)^{-1}\!.\]
\end{enumerate}
\end{theorem}

In the simplest flopping situation, where $f$ contracts a $(-1,-1)$-curve, the autoequivalence $\twistLocal$ is a spherical twist in the sense of Seidel--Thomas~\cite{ST}. For a contraction of a $(-2,0)$-curve, it is a generalized spherical twist as first constructed by Toda~\cite{Toda}, who furthermore established the conclusion of Theorem~\ref{theorem flop flop}(\ref{theorem flop flop 3}) in this case.

\begin{remark} The noncommutative deformation theory used here relies on work of Laudal~\cite{Laudal}, Eriksen~\cite{Eriksen}, E.\ Segal~\cite{Segal}, and Efimov--Lunts--Orlov~\cite{ELO1}.
\end{remark}

\begin{remark} The algebra $A$ above, and similar noncommutative deformation algebras, have now been applied in settings including: enumerative geometry of curves on $3$-folds by Toda and Hua--Toda~\cite{TodaGV,HuaT}; flops of families of curves in higher dimensions by Bodzenta and Bondal~\cite{BB};  construction of  autoequivalences and exceptional objects by Kawamata~\cite{Kawamata}; and  new braid-type groups of derived symmetries of $3$-folds by the author and Wemyss~\cite{DW3}.
 \end{remark}

\begin{remark} The full statement of Theorem~\ref{theorem flop flop} does not require $X$ to be smooth: I leave details to the references.
\end{remark}

\section{General results} 

The following theorem from~\cite{DW4} gives a sheafy analogue of the deformation algebra $A$, applicable in higher dimensions. For a birational contraction $f\colon X \to Y$ satisfying the assumption below, we define a sheaf of algebras $\mathcal{A}$ on $Y$ which is supported on the locus over which $f$ is not an isomorphism. We furthermore construct an associated autoequivalence of $D(X)$.

\begin{assumption} Suppose that $f\colon X \to Y$ is a contraction with $\operatorname{dim} X \geq 2$, and that either:
\begin{enumerate}[label=(\alph*),ref=\alph*]
\item\label{assumption 1} the variety $X$ has an $f$\!-relative tilting generator with summand $\mathcal{O}_X$, where $f$ is crepant, and \,$Y$\!~is Gorenstein;
\end{enumerate}
or, alternatively,
\begin{enumerate}[label=(\alph*),ref=\alph*,resume]
\item\label{assumption 2} the fibres of $f$ have dimension at most one.
\end{enumerate}
\end{assumption}

\begin{remark} The tilting generator assumption from~(\ref{assumption 1}) is satisfied in a range of situations, including symplectic resolutions of quotient singularities as established by Bezrukavnikov and Kaledin~\cite{BK}, and contractions with fibres of dimension at most two under conditions of Toda and~Uehara~\cite{TodaUehara}. 
\end{remark}

Write $Z$ for the locus in $Y$\! over which $f$ is not an isomorphism. 

\begin{figure}[h]
\begin{center}
\begin{tikzpicture}
\node at (0,0) {\begin{tikzpicture}[scale=0.6]
\coordinate (T) at (1.9,2);
\coordinate (B) at (2.1,1);

\coordinate (sT) at (1.9,2.3);
\coordinate (sL) at (1.7,1.5);
\coordinate (sR) at (2.35,1.5);
\coordinate (sB) at (2.2,0.7);

\draw[rounded corners=7pt,line width=0.5pt] (sT) -- (sR) -- (sB) -- (sL) -- cycle;
\foreach \x/\y/\h in {0.31/0.02/0.08,0.51/0.03/0.06,0.71/0.04/0.04,0.91/0.05/0.02}{ 
\draw[line width=0.5pt] ($(T)+(\x,\y)+(0.02,\h)$) to [bend left=25] ($(B)+(\x,\y)+(0.02,-\h)$);
\draw[line width=0.5pt] ($(T)+(-\x,-\y)+(0.02,\h)$) to [bend right=25] ($(B)+(-\x,-\y)+(0.02,-\h)$);
}
\draw[dashed,rounded corners=6pt,line width=0.75pt] (0.5,0.3) -- (1.5,0.3)-- (3.6,0) -- (4.3,1.5)-- (3.8,2.7) -- (2.5,2.7) -- (0.3,2.8) -- (-0.2,1.8)-- cycle;
\node at (-0.5,2.5) {$X$};
\end{tikzpicture}};
\node at (0,-2) {\begin{tikzpicture}[scale=0.6]
\filldraw (2.1,0.81) circle (1pt);
\draw (1,0.65) to [bend left=5] (3,0.85);
\draw[dashed,rounded corners=5pt,line width=0.75pt] (0.5,0.15) -- (1.5,0.15)-- (3.6,0) -- (4.3,0.75)-- (3.8,1.35) -- (2.5,1.35) -- (0.3,1.4) -- (-0.2,0.9)-- cycle;
\node at (-0.5,0.5) {$Y$};
\node (label) at (3.4,0.85) {$Z$};
\end{tikzpicture}};
\draw[->] (0.25,-0.9) -- node[color=black,left]{$f$} (0.25,-1.4);
\end{tikzpicture}
\end{center}
\caption{Contraction $f$ for Theorem~\ref{theorem dw4 simple version}.}
\end{figure}
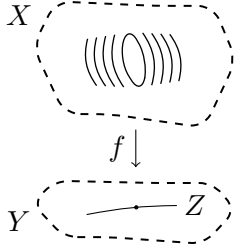

\begin{theorem}\label{theorem dw4 simple version} \cite{DW4}  Under the assumption above, there is a sheaf of algebras~$\mathcal{A}$ on~$Y$\!, inducing an object $\mathcal{E}$ of $D(X)$, such that:
\begin{enumerate}
\item\label{theorem dw4 simple version 1} the support of $\mathcal{A}$ is $Z$.
\end{enumerate}
For points $z$ of $Z$ such that $f^{-1}(z)$ is one-dimensional with components~$C_i$, then:
\begin{enumerate}[start=2]
\item\label{theorem dw4 simple version 3} the completion $\mathcal{A}_z$ is an algebra which prorepresents the functor of noncommutative deformations of the sheaves $\mathcal{O}_{C_i}(-1)$ on $X$, up to Morita equivalence;
\item\label{theorem dw4 simple version 4} the restriction of\, $\mathcal{E}$ to the formal fibre of $f$\! over $z$ is a sheaf, namely the universal family corresponding to the prorepresenting object given in~(\ref{theorem dw4 simple version 3}), up to summands of finite sums of sheaves.
\end{enumerate}
If the following hold:
\begin{enumerate}[label=(\roman*),ref=\roman*]
\item\label{theorem dw4 simple version assumptions 1} the contraction $f$ is crepant,
\item\label{theorem dw4 simple version assumptions 2} the base $Y$ is complete locally a hypersurface at each point of~$Z$, 
\end{enumerate} 
and either $\operatorname{codim} Z \geq 3$ or, alternatively,
\begin{enumerate}[label=(\roman*),ref=\roman*,resume]
\item\label{theorem dw4 simple version assumptions 3} the sheaf $\mathcal{A}$ is Cohen--Macaulay, and
\item\label{theorem dw4 simple version assumptions 4} the object $\mathcal{E}$ is perfect,
\end{enumerate} 
then:
\begin{enumerate}[start=4]
\item\label{theorem dw4 simple version 2} there is a Fourier--Mukai autoequivalence~$\twistGlobal$ of~$D(X)$, fitting into a distinguished triangle of functors
\begin{equation*}\label{triangle 2} f^{-1} \mathbb{R}f_* \mathbb{R}\mathcal{H}om_X(\mathcal{E},-) \!\!\underset{f^{-1} \mathcal{A}}{\overset{\mathbb{L}}{\otimes}}\! \mathcal{E} \longrightarrow \operatorname{Id}_{D(X)} \longrightarrow \twistGlobal \longrightarrow.
\end{equation*}
\end{enumerate}
\end{theorem}

I indicate the construction of the sheaf of algebras $\mathcal{A}$, and explain how it allows us to prove Theorem~\ref{theorem dw4 simple version}. Under the assumption above, we have an $f$-relative tilting generator $\mathcal{O}_X \oplus N$, either by assertion in case~(\ref{assumption 1}), or by a theorem of Van~den~Bergh~\cite{VdB} in case~(\ref{assumption 2}). Let $\mathcal{T}$ denote the relative endomorphism algebra of $\mathcal{O}_X \oplus N$, a~sheaf of algebras on~$Y$. We establish that
\[\mathcal{T} =   f_* \mathcal{E}nd_X (\mathcal{O}_X \oplus N) \cong \mathcal{E}nd_Y f_* (\mathcal{O}_X \oplus N). \]
This allows us to make the following definition for $\mathcal{A}$. This is a sheafy version of a construction of the algebra $A$ from our previous work~\cite{DW1}.
\begin{definition}\cite{DW4} Let $\mathcal{A} = \mathcal{T} / \mathcal{I}$, a sheaf of algebras on $Y$\!, where $\mathcal{I}$ is the ideal of sections of $\mathcal{T}$ which factor, at each stalk, through a sum of copies of $\mathcal{O}_Y$.
\end{definition}
The prorepresenting property of $\mathcal{A}$ in Theorem~\ref{theorem dw4 simple version}(\ref{theorem dw4 simple version 3})  is then proved as a sheafy version of the representing property of $A$ in Theorem~\ref{theorem flop flop}(\ref{theorem flop flop 1}). The object $\mathcal{E}$ of $D(X)$ is defined as the image of $\mathcal{A}$ under an appropriate tilting equivalence: I refer to \cite[Section~3]{DW4} for a precise statement. Theorem~\ref{theorem dw4 simple version}(\ref{theorem dw4 simple version 4}) show that this $\mathcal{E}$ has a universal property which is a sheafy version of the universal property of $E$ in Theorem~\ref{theorem flop flop}. We also have the following.

\begin{proposition}\cite{DW4} The support of $\mathcal{E}$ is contained in the exceptional locus of~$f$.
\end{proposition}

The construction of a Fourier-Mukai autoequivalence~$\twistGlobal$ in Theorem~\ref{theorem dw4 simple version}(\ref{theorem dw4 simple version 2}) generalizes the construction of~$\twistLocal$ from Theorem~\ref{theorem flop flop}(\ref{theorem flop flop 2}). In particular, we have the following.

\begin{remark}When $Z$ is a point, the autoequivalence~$\twistGlobal$  reduces to the autoequivalence~$\twistLocal$ appearing in Theorem~\ref{theorem flop flop}(\ref{theorem flop flop 2}).
\end{remark}

\begin{remark}{} Although the tilting generator $\mathcal{O}_X \oplus N$, and thence the sheaf of algebras $\mathcal{A}$, is not canonically defined (see for instance the construction of Van~den~Bergh in~\cite{VdB}) it seems that the autoequivalence $\twistGlobal$ may be canonical, given a choice of contraction~$f$. For instance, given two different tilting generators related by duplication of summands we obtain Morita equivalent sheaves of algebras $\mathcal{A}$, and thence isomorphic autoequivalences~$\twistGlobal$. \end{remark}

\begin{remark}
For $f$ a flopping contraction, it would be interesting to establish when $\twistGlobal$ is related to a flop-flop functor, as in Theorem~\ref{theorem flop flop}(\ref{theorem flop flop 3}).
\end{remark}

\begin{remark}
It is tempting to speculate that the `tilting' condition in the requirement for an $f$-relative tilting generator in~(\ref{assumption 1}) may be relaxed by upgrading~$\mathcal{A}$ to an appropriate sheaf of differential graded algebras.
\end{remark}

I record the following $3$-fold setting where the assumptions of Theorem~\ref{theorem dw4 simple version} may be established.

\begin{theorem}\label{theorem dw4 3fold irred version}\cite{DW4} With $\operatorname{dim} X = 3$,  assume that

\begin{enumerate}[label=(\roman*),ref=\roman*]
\item the contraction $f$ is crepant,
\item the base $Y$ is complete locally a hypersurface at each point of~$Z$,
\item the exceptional fibres of $f$ are irreducible curves.  
\end{enumerate} 
Then the assumptions of Theorem\! \ref{theorem dw4 simple version} hold, and  there exists an associated autoequivalence $\twistGlobal$ of $D(X)$.
\end{theorem}

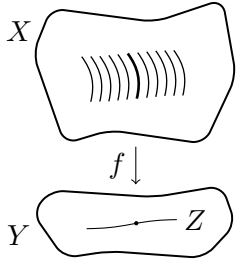
\begin{figure}[h]
\begin{center}
\begin{tikzpicture}
\node at (0,0) {\begin{tikzpicture}[scale=0.6]
\coordinate (T) at (1.9,2);
\coordinate (TM) at (2.12-0.02,1.5-0.15);
\coordinate (BM) at (2.12-0.09,1.5+0.2);
\coordinate (B) at (2.1,1);
\draw[line width=1pt] (T) to [bend left=25] (B);
\foreach \x/\y in {0.2/0.02,0.4/0.04,0.6/0.06,0.8/0.065,1/0.07}{ 
\draw[line width=0.5pt] ($(T)+(\x,0)+(0.02,\y)$) to [bend left=25] ($(B)+(\x,0)+(0.02,\y)$);
\draw[line width=0.5pt] ($(T)+(-\x,0)+(-0.02,-\y)$) to [bend left=25] ($(B)+(-\x,0)+(-0.02,-\y)$);}
\draw[rounded corners=5pt,line width=0.75pt] (0.5,0) -- (1.5,0.3)-- (3.6,0) -- (4.3,1.5)-- (4,3.2) -- (2.5,2.7) -- (0.2,3) -- (-0.2,2)-- cycle;
\node at (-0.5,2.5) {$X$};
\end{tikzpicture}};
\node at (0,-2) {\begin{tikzpicture}[scale=0.6]
\filldraw (2.1,0.77) circle (1pt);
\draw (1,0.65) to [bend right=5] (2,0.75);
\draw (2,0.75) to [bend left=5] (3,0.85);
\draw[rounded corners=5pt,line width=0.75pt] (0.5,0) -- (1.5,0.15)-- (3.6,0) -- (4.3,0.75)-- (4,1.6) -- (2.5,1.35) -- (0.2,1.5) -- (-0.2,1)-- cycle;
\node at (-0.5,0.5) {$Y$};
\node (label) at (3.4,0.85) {$Z$};
\end{tikzpicture}};
\draw[->] (0.25,-1) -- node[color=black,left]{$f$} (0.25,-1.5);
\end{tikzpicture}
\end{center}
\caption{Contraction $f$ for Theorem~\ref{theorem dw4 3fold irred version}.}
\end{figure}

\newpage
\section{Springer resolution example}\label{section springer}

For an example in which the theory of the previous section applies to a contraction with higher-dimensional fibres, consider the Springer resolution of the variety of singular $d$-by-$d$ matrices. Namely, for a vector~space~$V$ of dimension~$d$ with $d \geq 2$, take the singular cone
\[ Y = \big\{ M \in \operatorname{End} V  \mathrel{\big|} \det M = 0 \,\big\}, \]
which is a Gorenstein hypersurface. It has a resolution by 
\[ X = \big\{ (M,H) \in \operatorname{End} V \times \mathbb{P} V \mathrel{\big|} \operatorname{Im} M \subseteq H \,\big\} \] 
whose natural projection\! $f$ to $\operatorname{End} V$ surjects onto $Y$. This resolution~$f$ is crepant. Its exceptional fibres lie over points $M$ in $Y$ with $\operatorname{rk} M < d-1$, and are projective spaces of dimension $d-1-\operatorname{rk} M$.

A tilting generator for $X$ has been constructed by Buchweitz, Leuschke, and Van~den~Bergh \cite{BLV}, so that we are in the setting of Assumption~(\ref{assumption 1}). Conditions~(\ref{theorem dw4 simple version assumptions 1}) and~(\ref{theorem dw4 simple version assumptions 2}) of Theorem~\ref{theorem dw4 simple version} are noted above, and 
$\operatorname{codim} Z \geq 3$, so we can apply Theorem~\ref{theorem dw4 simple version}(\ref{theorem dw4 simple version 2}) to obtain an autoequivalence $\twistGlobal$ of $D(X)$.

\begin{remark}For $d=2$, the variety $X$ is just a $3$-fold resolving a conifold~$Y$ with exceptional fibre a $(-1,-1)$-curve, and we are in the setting of Theorem~\ref{theorem flop flop}.
\end{remark}

\begin{remark}
In joint work with E.~Segal~\cite{DonSeg1}, I studied the resolution $f$, along with more general resolutions where $Y$ is the variety of $d$-by-$d$ matrices of rank at most~$r$ for $0 < r < d$. We constructed certain `Grassmannian twist' autoequivalences of the corresponding $D(X)$ by quite different methods: it~would be interesting to compare these with $\twistGlobal$. 
\end{remark}

\begin{remark} The sheaf of algebras $\mathcal{A}$ for this example may be computed from the presentation of the endomorphism algebra $\mathcal{T}$ in~\cite{BLV}.
\end{remark}

\section{Conjectures}\label{section conjectures}

In the setting of a $3$-fold flopping contraction as in Theorem~\ref{theorem flop flop}, we made a conjecture~\cite[Conjecture~1.4]{DW1} that the complete local neighbourhood of the $3$-fold $Y$ near the singularity~$p$ is determined, up to isomorphism, by~the deformation algebra $A$. This conjecture is clear in the following simple cases, namely the two kinds of flopping curve for which $A$ is commutative, but remains open more generally.

\newpage\begin{enumerate}
\item Contractions of $(-1,-1)$-curves. In this case $A=\mathbb{C}$, and the completion of $Y$ at $p$ is necessarily a conifold singularity.
\item Contractions of $(-2,0)$-curves. Here $A \cong \mathbb{C}[\varepsilon]/\varepsilon^w$ with $w\geq 2$, where~$w$ is the width invariant of Reid~\cite{Pagoda}, and the completion of~$Y$ at $p$ is determined by $w$.
\end{enumerate}

Hua and Toda subsequently proposed an $A_\infty$ version of the conjecture~\cite[Conjecture~5.3]{HuaT} in which $A$ is endowed with the structure of an $A_\infty$-algebra. They established their conjecture for contractions to weighted homogeneous hypersurface singularities~\cite[Theorem~5.5]{HuaT}, and it has now been settled in general by Hua~\cite{Hua}. A key idea in these works is that the deformation algebra~$A$ may be viewed as a noncommutative analogue of the Milnor algebra, and that the $A_\infty$ structure on it allows recovery of the Milnor algebra along with enough information to apply a version of the Mather--Yau theorem.

\begin{remark}
It would be interesting to extend these results to higher dimensions, and to non-isolated singularities. For instance, it is natural to ask whether the the complete local neighbourhood of the variety $Y$\! near the non-isomorphism locus $Z$ is determined by $(Z,\mathcal{A})$, potentially along with some appropriate $A_\infty$~structure.\end{remark}

\newpage
\renewcommand{\baselinestretch}{1.15}


\begin{thebibliography}{99}
\bibitem{BK}
R.~Bezrukavnikov and D.~Kaledin, \emph{McKay equivalence for symplectic resolutions of singularities}, Tr. Mat. Inst. Steklova 246 (2004), Algebr. Geom. Metody, Svyazi i Prilozh., 20--42; English transl. Proc Steklov Inst. Math. 246 (2004), 13--33.

\bibitem{BB}
A.~Bodzenta and A.~Bondal, \emph{Flops and spherical functors}, \href{http://arxiv.org/abs/1511.00665}{\texttt{arXiv:1511.00665}}.

\bibitem{BriFlop}
T.~Bridgeland, \emph{Flops and derived categories}, Invent. Math. \textbf{147} (2002), no.~3, 613--632.

\bibitem{BLV}
R-O.~Buchweitz, G.~J.~Leuschke and M.~Van den Bergh, \emph{Non-commutative desingularization of determinantal varieties I}, Invent.\ Math.\ \textbf{182} (2010), no.~1, 47--115.

\bibitem{DonSeg1}
W.~Donovan and E.~Segal, \emph{Window shifts, flop equivalences and Grassmannian twists}, Compos.\ Math.\ (6) \textbf{150} (2014), 942--978.

\bibitem{DW1}
W.~Donovan and M.~Wemyss, \emph{Noncommutative deformations and flops}, Duke Math.~J. \textbf{165} (8) (2016): 1397--1474, \href{http://arxiv.org/abs/1309.0698}{\texttt{arXiv:1309.0698}}.

\bibitem{DW2}
\bysame, \emph{Contractions and deformations}, \href{http://arxiv.org/abs/1511.00406}{\texttt{arXiv:1511.00406}}.

\bibitem{DW3}
\bysame, \emph{Twists and braids for general 3-fold flops}, to appear {{J.\ Eur.\ Math.\ Soc.}}, \href{http://arxiv.org/abs/1504.05320}{\texttt{arXiv:1504.05320}}.

\bibitem{DW4}
\bysame, \emph{Noncommutative enhancements of contractions}, \href{http://arxiv.org/abs/1612.01687}{\texttt{arXiv:1612.01687}}.

\bibitem{ELO1}
A.~Efimov, V.~Lunts and D.~Orlov, \emph{Deformation theory of objects in homotopy and derived categories I: General theory}, Adv.\ Math.\ \textbf{222} (2009), no.~2, 359--401.


\bibitem{Eriksen}
E.~Eriksen, \emph{Computing noncommutative deformations of presheaves and sheaves of modules}, Canad.\ J.\ Math.\ \textbf{62} (2010), no.~3, 520--542.

\bibitem{Hua}
Z.~Hua, \emph{On a conjecture of Donovan and Wemyss}, \href{http://arxiv.org/abs/1610.05467}{\texttt{arXiv:1610.05467}}.

\bibitem{HuaT}
Z.~Hua and Y.~Toda, \emph{Contraction algebra and invariants of singularities}, \href{http://arxiv.org/abs/1601.04881}{\texttt{arXiv:1601.04881}}.

\bibitem{Kawamata}
Y.~Kawamata, \emph{On multi-pointed non-commutative deformations and Calabi-Yau threefolds}, \href{http://arxiv.org/abs/1512.06170}{\texttt{arXiv:1512.06170}}.

\bibitem{Laudal}
O.~A.~Laudal, \emph{Noncommutative deformations of modules}. The Roos Festschrift volume, 2.\ Homology Homotopy Appl.\ \textbf{4} (2002), no.~2, part 2, 357--396.

\bibitem{Pagoda}
M.~Reid, \emph{Minimal models of canonical 3-folds}, Algebraic varieties and analytic varieties (Tokyo, 1981), 131--180, Adv.\ Stud.\ Pure Math., 1, North-Holland, Amsterdam, 1983.

\bibitem{Segal}
E.~Segal, \emph{The $A_\infty$ deformation theory of a point and the derived categories of local Calabi--Yaus}, J.\ Alg.\ \textbf{320} (2008), no.~8, 3232--3268.


\bibitem{ST}
P.~Seidel and R.~P.~Thomas, \emph{Braid group actions on derived categories of sheaves}, Duke Math.\ J.\ \textbf{108} (2001), 37--108.

\bibitem{Toda}
Y.~Toda, \emph{On a certain generalization of spherical twists},
Bulletin de la Soci\'et\'e Math\'ematique de France \textbf{135}, fascicule 1 (2007), 119--134.

\bibitem{TodaGV} 
\bysame, \emph{Non-commutative width and Gopakumar--Vafa invariants}, Manuscripta Math.\ \textbf{148} (2015), no.~3--4, 521--533.

\bibitem{TodaUehara}
Y.~Toda and H.~Uehara, \emph{Tilting generators via ample line bundles}, Adv.\ Math.\ \textbf{223} (2010), no.~1, 1--29.

\bibitem{VdB}
M.~Van den Bergh, \emph{Three-dimensional flops and noncommutative rings}, 
Duke Math.~J.\ \textbf{122} (2004), no.~3, 423--455. 

\end{thebibliography}
\end{document}